\documentclass[11pt]{article}
\usepackage{epic,latexsym,amssymb}
\usepackage{amsfonts}
\usepackage{amscd}
\usepackage{amsmath}
\usepackage{graphicx}
\usepackage{color}
\usepackage{caption,subcaption}
\usepackage{tikz}

\usepackage{float}

\newtheorem{thm}{Theorem}
\newtheorem{conj}{Conjecture}
\newtheorem{ob}[thm]{Observation}
\newtheorem{prop}[thm]{Proposition}

\newtheorem{quest}{Question}

\newcommand{\diam}{{\rm diam}}

\newcommand{\cP}{\mathcal{P}}

\newcommand{\barD}{\overline{D}}

\newcommand{\proof}{\noindent\textbf{Proof. }}

\newcommand{\smallqed}{{\tiny ($\Box$)}}
\newcommand{\qed}{$\Box$}

\newcommand{\2}{\vspace{0.2cm}}

\newenvironment{unnumbered}[1]{\trivlist
\item [\hskip \labelsep {\bf #1}]\ignorespaces\it}{\endtrivlist}

\newcommand{\QEDmark}{\mbox{\textsc{qed}}}
\newcommand{\proofStarter}[1]{\textsc{#1} }

\def\vertex(#1){\put(#1){\circle*{2}}}
\def\vertexo(#1){\put(#1){\circle{2}}}
\def\vert(#1){\put(#1){\circle*{1.5}}}
\def\verto(#1){\put(#1){\circle{1.5}}}
\def\lab(#1)#2{\put(#1){\makebox(0,0)[c]{#2}}}
\definecolor{DarkGreen}{rgb}{0.2, 0.6, 0.3}

\definecolor{electricindigo}{rgb}{0.44, 0.0, 1.0}

\begin{document}

\title{Bounds on the Connected Forcing Number \\ of a Graph}
\author{$^{1,4}$Randy Davila, $^1$Michael A. Henning, $^2$Colton Magnant, and $^3$Ryan Pepper \\
\\
$^1$Department of Pure and Applied Mathematics\\
University of Johannesburg \\
Auckland Park 2006, South Africa \\
\small {\tt Email: mahenning@uj.ac.za} \\
\\
$^2$Department of Mathematical Sciences \\
Georgia Southern University \\
Statesboro, GA 30458, USA \\
\small {\tt Email: cmagnant@georgiasouthern.edu}\\
\\
$^3$Department of Mathematics \\
University of Houston-Downtown \\
Houston, TX 77002, USA \\
\small {\tt Email: pepperr@uhd.edu}\\
\\
$^4$Department of Mathematics \\
Texas State University \\
San Marcos, TX 78666, USA \\
\small {\tt Email: rrd32@txstate.edu}
}

\date{}
\maketitle

\begin{abstract}
In this paper, we study (zero) forcing sets which induce connected subgraphs of a graph.  The minimum cardinality of such a set is called the connected forcing number of the graph.  We provide sharp upper and lower bounds on the connected forcing number in terms of the minimum degree, maximum degree, girth, and order of the graph.
\end{abstract}


\section{Introduction}
Graph coloring is one of the most widely applied and studied concepts in graph theory, but for the duration of the 20th century almost all graph colorings were static. That is, a given graph coloring would not change over the course of time. However, in recent years, new variations of graph colorings have emerged that allow a coloring of a graph to adapt over time in discrete intervals. A graph may then acquire many different colorings based on a single initial coloring. One of the most prominent dynamic coloring processes is called the \emph{forcing process} (\emph{zero forcing process}) with its associated graph parameter, the \emph{forcing number} (\emph{zero forcing number}). These concepts were originally introduced at a workshop on linear algebra and graph theory in 2006 \cite{AIM-Workshop} and quickly found a variety of applications in physics, logic circuits, coding theory, and power network monitoring \cite{quantum1, logic1, powerdom3, powerdom2}.

Let $G=(V,E)$ be a connected simple graph with vertex set $V=V(G)$ and edge set $E=E(G)$. The forcing process is defined as follows: Let $S \subseteq V$ be an initial set of ``colored'' vertices; all remaining vertices being ``non-colored''. A vertex in a set $S$, we call $S$-\emph{colored}, while a vertex not in $S$ we call $S$-\emph{uncolored}. At each time step, a colored vertex with exactly one non-colored neighbor will change, or {\it force}, the non-colored neighbor to be colored. We call such a vertex a \emph{forcing colored vertex}, or simply a \emph{forcing vertex}. A set $S \subseteq V$ of initially colored vertices is called a \emph{forcing set} if, by iteratively applying the forcing process, all of $V$ becomes colored. We call such a set $S$ an  $S$-\emph{forcing set}. The \emph{forcing number} of a graph $G$, denoted by $F(G)$, is the cardinality of a smallest forcing set. If $S$ is a forcing set in $G$ and $v$ is an $S$-colored vertex that forces a new vertex to be colored, then we call $v$ an $S$-\emph{forcing vertex}.

In general, the problem of determining $F(G)$ is in the class of $NP$-hard problems \cite{ sdmr_tf,zf_np}. If $S$ is a forcing set that induces a connected subgraph, we say that $S$ is a \emph{connected forcing set}. The cardinality of a smallest connected forcing set in $G$ is its \emph{connected forcing number}, denoted $F_c(G)$, and provides a new graph invariant that we introduce \footnote{This invariant is concurrently introduced in \cite{Brimkov Davila}.} and study in this paper.

We denote the order and size of $G$, by $n=n(G)=|V(G)|$ and $m=m(G)=|E(G)|$, respectively. Two vertices $v,w\in V$ are said to be adjacent, or neighbors, if there exists the edge $vw\in E$. The open neighborhood of $v\in V$ is the set of all vertices which are adjacent to $v$, denoted  $N(v)$. The closed neighborhood of $v\in V$, is $N[v]=N(v)\cup\{v\}$. Similarly, we define the open and closed neighborhoods of $S\subseteq V$, to be $N(S):=\{w: w\in N(v) \:\text{and} \: v\in S\}\setminus S$ and $N[S]=\{w: w\in N[v] \:\text{and} \: v\in S\}$, respectively. The degree of $v\in V$ is defined as $d(v)=|N(v)|$. The minimum degree and maximum degree of $G$ will be denoted as $\delta=\delta(G)$ and $\Delta=\Delta(G)$, respectively.

Let $A$ and $B$ be vertex disjoint subsets of vertices in a graph $G$. The set of edges between $A$ and $B$ in $G$ is denoted by $[A,B]$. If $S$ is a subset of vertices in $G$ and if $v \in V(G)$, then the \emph{degree of $v$ in $S$}, denoted by $d_S(v)$, is the number of vertices in $S$ adjacent to $v$. In particular, if $S = V(G)$, then $d_S(v)$ is the degree, $d_G(v)$, of $v$ in $G$.

For two vertices $u$ and $v$ in a connected graph $G$, the \emph{distance} $d_G(u,v)$ between $u$ and $v$ is the length of a shortest $(u,v)$-path in $G$. The maximum distance among the vertices of $G$ is its \emph{diameter}, which is denoted by $\diam(G)$. For a set $S \subseteq V$, we let $G[S]$ denote the subgraph induced by $S$. The length of a shortest cycle in $G$ is the girth of $G$, denoted $g=g(G)$. We denote a path and cycle on $n$ vertices by $P_n$ and $C_n$, respectively.

\paragraph{Domination in graphs.}  A set $S\subseteq V$ of vertices in a graph $G$ is a \emph{dominating set} if every
vertex not in $S$ is adjacent to some vertex in $S$. If a dominating set induces a connected subgraph, then we say that it is a \emph{connected dominating set}. The \emph{domination number} is the cardinality of a minimum dominating set, denoted $\gamma(G)$. The \emph{connected domination} number is the cardinality of a minimum connected dominating set, denoted $\gamma_c(G)$. As another variant of domination, the power domination process is defined as follows. For a given set $S\subseteq V$, the sets $\big(\cP_{G}^{i}(S)\big)_{i\ge 0}$ of vertices \emph{monitored} by $S$ at the $i$-th step are defined recursively by,
\begin{enumerate}
\item $\cP_{G}^{0}=N[S]$, and
\item $\cP_{G}^{i+1} = \bigcup \{ N[v]\: : \: v\in \cP_{G}^{i}(S)\:\ \text{such that} \:\: |N[v]\setminus \cP_{G}^{i}(S)|\le 1\}$.
\end{enumerate}
If $\cP_{G}^{i_{0}}=\cP_{G}^{i_{0}+1}$, for some $i_{0}$, then $\cP_{G}^{j}=\cP_{G}^{i_{0}}$, for all $j\ge i_{0}$. We define $\cP_{G}^{\infty}=\cP_{G}^{i_{0}}$. If $\cP_{G}^{\infty}(S)=V$, we say that $S$ is a \emph{power dominating set} of $G$. The cardinality of a smallest power dominating set is the \emph{power domination number} of $G$, and is denoted $\gamma_P(G)$.

\paragraph{Known Results on Forcing Domination.} It should be highlighted that $F(G)$ is, in general, very difficult to compute, even for well structured graphs like bipartite graphs. The difficulty in computing $F(G)$ has motivated mathematicians to seek computationally efficient upper and lower bounds on $F(G)$.

In a paper on generalized forcing,  Amos, Caro, Davila, and Pepper~\cite{k-Forcing} established an upper bound on the forcing number of a graph in terms of the order $n$ and the maximum degree~$\Delta$. In particular, they proved $F(G) \le (\frac{\Delta}{\Delta +1})n$ for isolate-free graphs $G$, and $F(G) \leq\frac{(\Delta-2)n+2}{\Delta-1}$ for connected graphs $G$ with $\Delta(G) \ge 2$. These results resolved a question posed by Meyer \cite{Meyer} which asked if the forcing number could be bounded from above by a function of order and degree for bipartite circulant graphs. Moreover, Amos et al.~\cite{k-Forcing} also related the forcing number to the connected domination number with the inequality $F(G)\le n - \gamma_c(G)$. A slight improvement on the work of Amos et al., Caro and Pepper \cite{Dynamic Forcing} used a greedy algorithm on connected graphs with maximum degree $\Delta \ge 2$ to show $F(G)\le \frac{(\Delta -2)n-(\Delta - \delta) + 2}{\Delta -1}$.

The first lower bound in terms of the minimum degree came from the original paper on forcing due to the AIM-Group  \cite{AIM-Workshop}, which showed $F(G) \ge \delta$ for all graphs $G$. Improving upon this minimum degree lower bound, Davila and Kenter~\cite{Davila Kenter} proved $F(G) \ge 2\delta -2$ for graphs $G$ with girth $g \ge 5$, and $F(G) \ge \delta + 1$ for triangle-free graphs~$G$ with minimum degree $\delta \ge  3$. Further, Davila and Kenter conjectured $F(G) \ge \delta + (\delta-2)(g-3)$ which remains open.

We have two immediate aims in this paper. Our first aim is to introduce the connected forcing number of a graph and establish fundamental properties of this parameter. Our second aim is to provide lower and upper bounds on the connected forcing number, similar to those shown for the forcing number.

\section{Preliminary Observations and Results}

Since every isolated vertex in a graph must belong to every forcing set in the graph, we assume throughout this paper that there are no isolates. Since every connected forcing set is also a forcing set, we have the following observation.

\begin{ob}
\label{ob:trivial1}
For every connected graph $G$, it holds that $F(G) \le F_c(G)$.
\end{ob}

We first determine the connected forcing number of simple classes of graphs. In a complete graph on $n \ge 2$ vertices, no set of $n-2$ vertices is a forcing set, implying that $F(K_n) = F_c(K_n) = n-1$. Since the leaf of every non-trivial path is a connected forcing set, we note that $F(P_n) = F_c(P_n) = 1$. No vertex of a cycle is a forcing set, while any two adjacent vertices in a cycle form a forcing set in the cycle, implying that for $n \ge 3$, $F(C_n) = F_c(C_n) = 2$. We state these results formally as follows.

\begin{ob}
\label{o:PnCnKn}
For $n \ge 3$, the following hold. \\
\indent {\rm (a)} $F_c(P_n) = 1$. \\
\indent {\rm (b)} $F_c(C_n) = 2$. \\
\indent {\rm (c)} $F_c(K_n) = n-1$.
\end{ob}

If $G$ is a connected graph of order~$n \ge 2$, and if $v$ is a vertex of minimum degree in $G$, then $S \setminus \{v\}$ is a connected forcing set in $G$, implying that $F_c(G) \le |S| - 1 = n - 1$.

\begin{ob}
\label{o:trivial2}
If $G$ is a connected graph of order~$n \ge 2$, then $F_c(G) \le n - 1$.
\end{ob}

Suppose that $G$ is a connected graph of order~$n \ge 2$ satisfying $F(G) = 1$. Let $S$ be a minimum forcing set of $G$, and let $S = \{v_1\}$. Since $v_1$ is a forcing vertex, we note that $v_1$ is a leaf in $G$ and forces its neighbor, say $v_2$, to be colored. If $n > 2$, then the vertex $v_2$ has exactly one non-colored neighbor, implying that $v_2$ has degree~$2$ and forces its neighbor, say $v_3$, different from $v_1$ to be colored. If $n > 3$, then the vertex $v_3$ has exactly one non-colored neighbor, implying that $v_3$ has degree~$2$ and forces its neighbor, different from $v_2$, to be colored. Continuing in this way, the resulting graph $G$ is a path. We state this formally as follows.

\begin{ob}
\label{o:path}
Let $G$ be a connected graph of order~$n \ge 2$. Then, $F(G) = 1$ if and only if $G = P_n$.
\end{ob}

By Observation~\ref{ob:trivial1},~\ref{o:PnCnKn} and~\ref{o:path}, we note that if $G$ be a connected graph of order~$n \ge 2$ and $F_c(G) = 1$, then $G = P_n$.

\section{A Relation with Power Domination}

In this section, we relate the connected forcing number of a graph with its power domination number. The following proposition  relates the power domination process to the forcing process. In particular, we show that every power dominating set is a dominating set of a forcing set. This relation appeared in the thesis of Davila~\cite{Davila Thesis}.

\begin{prop}\label{Basic Relation}
Let $G$ be a graph. A subset of vertices $S\subseteq V$ is a power dominating set of $G$ if and only if $\cP_{G}^{1}(S)$ is a forcing set of $G$.
\end{prop}
\proof Let $S\subseteq V$ be a power dominating set of $G$. Then color $\cP_{G}^{1}(S)$, i.e., color $N[S]$. Then, either all of $V$ is colored, or there is some colored vertex $v$ such that $|N[v]\setminus N[S]| = 1$, i.e., $v$ has at exactly one non-colored neighbor and is a forcing vertex. This process will continue until we have reached a set equivalent to $\cP_{G}^{\infty}=V$, since $S$ was power dominating. Hence, $N[S]$ is a forcing set.

Conversely suppose $\cP_{G}^{1}(S)$ is a forcing set. Then either all of $V$ is colored, and $S$ is a dominating set, and hence also power dominating, or there is a vertex $v \in \cP_{G}^{1}(S)$, such that $v$ has exactly one non-colored neighbor, i.e., $|N[v]\setminus N[S]| = 1$. This is assured at each forcing step until all of $V$ is colored. Hence, $S$ must be power dominating.~\qed

\medskip
As an immediate consequence of Proposition~\ref{Basic Relation}, we have the following relation.

\begin{ob}
\label{Trivial Relation}
If $G$ is a connected graph, then $\gamma_P(G) \le F(G)\le F_c(G)$.
\end{ob}

A fundamental result in domination theory, is Ore's Theorem~\cite{Ore's Theorem} which states that the domination number of a graph without isolated vertices is at most one-half the order of the graph. Combining Ore's Theorem with Proposition~\ref{Basic Relation}, yields the following relation between the power domination and connected forcing numbers of a graph.

\begin{prop}
If $G$ is a connected graph of order at least~$3$ that is not a path, then $\gamma_P(G)\le \frac{1}{2}F_c(G)$,
and this bound is sharp.
\end{prop}
\proof Let $S\subset V$ be a minimum connected forcing set, and so $F_c(G) = |S|$. Since $G$ is not a path, Observation~\ref{o:path} implies that $F_c(G) \ge 2$. Thus, the connected forcing set $S$ has size at least~$2$ and therefore induces a (connected) subgraph without isolated vertices. Let $D$ be a minimum dominating set in the graph $G[S]$ induced by $S$. By Ore's Theorem, $|D| \le |S|/2$. Since $S$ is a forcing set of $G$ and $S \subseteq \cP_{G}^{1}(D)$, the set $\cP_{G}^{1}(D)$ is a forcing set of $G$, implying by  Proposition~\ref{Basic Relation} that the set $D$ is a power dominating set of $G$. Hence, $\gamma_P(G)\le |D| \le \frac{1}{2}|S| = \frac{1}{2}F_c(G)$. By Observation~\ref{o:PnCnKn}(b), every cycle has connected forcing number~$2$. Thus, since every cycle has power domination number~$1$, the bound is trivially sharp for cycles.~\qed


\section{Upper Bounds}

In this section we investigate upper bounds on the connected forcing number on a graph. We first determine the graphs that achieve equality in the upper bound of Observation~\ref{o:trivial2}.

\begin{thm}
\label{t:upperbd}
Let $G$ be a connected graph of order~$n \ge 2$. Then, $F_c(G) = n - 1$ if and only if $G$ is a complete graph, $K_n$, with $n \ge 2$ or a star, $K_{1,n-1}$, with $n \ge 4$.
\end{thm}
\proof If $G \cong K_n$ where $n \ge 2$ or $G \cong K_{1,n-1}$ where $n \ge 4$, then it is immediate that $F_c(G) = n - 1$. This establishes the sufficiency. To prove the necessity, suppose that $F_c(G) = n - 1$. We proceed by induction on $n \ge 2$. The base case when $n \in \{2,3\}$ is immediate. For the inductive hypothesis, let $n \ge 4$ and assume that if $G'$ is a connected graph of order~$n'$, where $2 \le n' < n$, satisfying  $F_c(G') = n' - 1$, then $G \cong K_{n'}$ or $G \cong K_{1,n'-1}$ where $n' \ge 4$. Let $G$ be a connected graph of order~$n$ satisfying $F_c(G) = n - 1$. If $\Delta(G) = 2$, then $G$ is a path or a cycle, and so, by Observation~\ref{o:PnCnKn}, $F_c(G) \le 2 < n-1$, a contradiction. Hence, $\Delta(G) \ge 3$. If $G$ is a complete graph or a star, then the desired result follows. Hence, we may assume that $G$ is neither a complete graph nor a star.

Suppose that $G$ is a tree. By assumption, $G$ is not a star, and so $\diam(G) \ge 3$. Let $u$ and $v$ be two vertices at maximum distance apart in $G$, and so $d_G(u,v) = \diam(G)$. Since $\diam(G) \ge 3$ and $G$ is a tree, the set $V(G) \setminus \{u,v\}$ is a connected forcing set in $G$, implying that $F_c(G) \le n - 2$, a contradiction. Hence, $G$ is not a tree, and therefore contains at least one vertex, $v$ say, of degree at least~$2$ that is not a cut vertex.
Let $G' = G - v$ have order~$n'$. By choice of the vertex~$v$, the graph $G'$ is connected. By Observation~\ref{o:trivial2}, $F_c(G') \le n' - 1$.

We show firstly that $F_c(G') \le n' - 2$. Suppose, to the contrary, that $F_c(G') = n' - 1$. By the inductive hypothesis, either $G \cong K_{n'}$ or $G \cong K_{1,n'-1}$ where $n' \ge 4$. Suppose that $G \cong K_{n'}$. By assumption, $G \ncong K_n$. Let $u$ be a neighbor of $v$ in $G$, and let $w$ be a vertex that is not a neighbor of $v$ in $G$. We now consider the set $S = V(G) \setminus \{u,v\}$. Since every neighbor of $w$, except for $u$, is colored, the vertex $w$ is a forcing vertex in the set $S$ that forces the vertex~$u$ to be colored. Once $u$ is colored, then $u$ becomes a forcing vertex in the resulting set $S \cup \{u\}$ and forces the vertex~$v$ to be colored. Thus, $S$ is a forcing set of $G$, implying that $F_c(G) \le n - 2$, a contradiction.

Suppose next that $G \cong K_{1,n'-1}$ where $n' \ge 4$. Let $z$ be the center of the star $G'$. Since $v$ has degree at least~$2$ in $G$, the vertex $v$ is adjacent to at least one leaf of $G'$, say $x$. Let $y$ be a leaf in $G'$ different from $x$, and consider the set $S = V(G) \setminus \{v,y\}$. Since $x$ has precisely two neighbors in $G$, namely $v$ and $z$, the colored vertex $x$ is a forcing vertex in the set $S$ that forces the vertex~$v$ to be colored. Once $v$ is colored, the vertex $z$ (possibly, $v$ and $z$ are neighbor) becomes a forcing vertex in the resulting set $S \cup \{v\}$ and forces the vertex~$y$ to be colored. Thus, $S$ is a forcing set of $G$, implying that $F_c(G) \le n - 2$, a contradiction. Since both cases produce a contradiction, we deduce that $F_c(G') \le n' - 2$.

Let $S'$ be a connected forcing set of the (connected) graph $G'$. By assumption, $|S| \le n' - 2$. Since each forcing vertex adds one new vertex to the set, and the connectivity of the resulting new set if preserved, there exists a set $T$ of vertices of $G'$ such that $T$ is a connected forcing set in $G'$, $S' \subseteq T$ and $|T| = n' - 2$. Let $a_1$ and $a_2$ be the two $T$-uncolored vertices in $G'$. Renaming $a_1$ and $a_2$, if necessary, we may assume that there is a vertex $t_1$ in the set $T$ that forces $a_1$ to be colored in $G'$, and that there is a vertex $t_2$ in the resulting set $T \cup \{a_1\}$ that forces $a_2$ to be colored in $G'$.

If $v$ is adjacent in $G$ to a vertex of $T$, then $T \cup \{v\}$ is a connected forcing set of $G$ of cardinality~$n-2$, implying that $F_c(G) \le n - 2$, a contradiction. Hence, $v$ has degree exactly~$2$ in $G$, and is adjacent in $G$ only to $a_1$ and $a_2$. We now consider the set $T$ in the graph $G$.

Since the vertex $a_1$ is the only $T$-uncolored neighbor of $t_1$ in $G'$, and since $t_1$ and $v$ are not neighbors in $G$, we note that the vertex $a_1$ is the only $T$-uncolored neighbor of $t_1$ in $G$, implying that $t_1$ is a $T$-forcing vertex in $G$ that forces the vertex $a_1$ to be colored.

If $t_2 \ne a_1$, then $a_2$ is the only $(T \cup \{a_1\})$-uncolored neighbor of $t_2$ in $G'$, implying that since $t_2$ and $v$ are not neighbors, $t_2$ is a forcing vertex in $G$ that forces the vertex $a_2$ to be colored. In the resulting set, $T \cup \{q_1,a_2\}$, both $a_1$ and $a_2$ are forcing vertices that force the vertex $v$ to be colored. Thus, in this case when $t_2 \ne a_1$, the set $T$ is a connected forcing set in $G$, implying that $F_c(G) \le n - 3$, a contradiction. Thus, $t_2 = a_1$.
If $a_2$ is adjacent to a vertex in $G'$ different from $a_1$, then we can choose $t_2 \ne a_1$, a contradiction. Hence, $a_1$ is the only neighbor of $a_2$ in $G'$. Thus, both $a_2$ and $v$ have degree~$2$ in $G$. Further, $N_G(a_2) = \{a_1,v\}$ and $N_G(a_1) = \{a_2,v\}$.

Let $v'$ be a vertex in $G'$ at maximum distance from~$a_1$. Since $G \ncong K_n$, we note that $v'$ is not a neighbor of $a_1$ in $G'$. We now consider the set $D = V(G) \setminus \{v,v'\}$. By our choice of the vertex $v'$, the graph $G[D]$ is connected. The vertex $a_2$ is a $D$-forcing vertex and forces the vertex $v$ to be colored. Further, any neighbor of $v'$ in $G$ is also a $D$-forcing vertex and forces the vertex $'$ to be colored. Thus, the set $D$ is a connected forcing set of $G$, implying that $F_c(G) \le n - 2$, a contradiction. This completes the proof of the theorem.~\qed

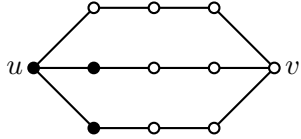
\begin{figure}[htb]

\begin{center}
\begin{tikzpicture}[scale=.8,style=thick,x=1cm,y=1cm]
\def\vr{2.5pt} 
\path (0,1) coordinate (v1);
\path (1,1) coordinate (v2);
\path (2,1) coordinate (v3);
\path (3,1) coordinate (v4);
\path (4,1) coordinate (v5);
\path (1,0) coordinate (u2);
\path (2,0) coordinate (u3);
\path (3,0) coordinate (u4);
\path (1,2) coordinate (w2);
\path (2,2) coordinate (w3);
\path (3,2) coordinate (w4);
\draw (v1) -- (u2);
\draw (v1) -- (v2);
\draw (v1) -- (w2);
\draw (v2) -- (v3);
\draw (v3) -- (v4);
\draw (v4) -- (v5);
\draw (u2) -- (u3);
\draw (u3) -- (u4);
\draw (w2) -- (w3);
\draw (w3) -- (w4);
\draw (v5) -- (u4);
\draw (v5) -- (w4);
\draw (u2) [fill=black] circle (\vr);
\draw (u3) [fill=white] circle (\vr);
\draw (u4) [fill=white] circle (\vr);
\draw (w2) [fill=white] circle (\vr);
\draw (w3) [fill=white] circle (\vr);
\draw (w4) [fill=white] circle (\vr);
\draw (v1) [fill=black] circle (\vr);
\draw (v2) [fill=black] circle (\vr);
\draw (v3) [fill=white] circle (\vr);
\draw (v4) [fill=white] circle (\vr);
\draw (v5) [fill=white] circle (\vr);
\draw[anchor = east] (v1) node {$u$};
\draw[anchor = west] (v5) node {$v$};
\end{tikzpicture}
\end{center}

\vskip -0.25 cm
\caption{A $3$-rail with shaded vertices colored} \label{Fig:Rail}
\end{figure}

\medskip
We establish next an upper bound on the connected forcing number of a $2$-connected graph in terms of its order and girth. For this purpose, we shall need the following definition. A $k$-\emph{rail} of a graph $G$ is a subgraph of $G$ consisting of two vertices with $k$ internally vertex disjoint paths between them such that the internal vertices of the paths have degree~$2$ not only in the $k$-rail but also in the graph $G$ itself. See Figure~\ref{Fig:Rail}.


\vskip 0.25cm

We shall need the following result due to Thomassen and Toft~\cite{TT81}.

\begin{thm}[\cite{TT81}] \label{Thm:TT}
If $G$ is a $2$-connected graph such that the removal of any induced cycle separates the graph into at least two components, then $G$ contains a $3$-rail.
\end{thm}

We are now in a position to prove the following result.

\begin{thm}
If $G$ is a $2$-connected graph of order $n$ with girth $g$, then $F_{c}(G) \le n - g + 2$.
\end{thm}
\proof
Among all induced cycles in $G$, let $C$ be chosen so that $G - V(C)$ contains as few components as possible. We note that $C$ has length at least~$g$. If $G - V(C)$ is connected, then color all of $V(G) \setminus V(C)$ along with two consecutive vertices of $C$ chosen to be adjacent to at least one vertex of $V(G) \setminus V(C)$. Let $S$ denote the resulting set of colored vertices. This colored set $S$ is a forcing set of $G$. Further, since $G - V(C)$ is connected, the set $S$ is a connected forcing set of $G$, implying that $F_{c}(G) \le |S| = n - |V(C)| + 2 \le n - g + 2$, as desired. Hence, we may assume that $G - V(C)$ is disconnected. Thus, by our choice of the cycle $C$, there is no induced cycle $C^*$ in $G$ such that $G - V(C^*)$ is connected.

By Theorem~\ref{Thm:TT}, there is a pair of vertices in $G$, say $\{u, v\}$, with a $3$-rail, say $R$, between them. Let $Q_1$, $Q_2$ and $Q_3$ be the three $(u,v)$-paths in $R$, where $Q_i$ has $q_i$ internal vertices for $i \in [3]$ and where $q_1 \le q_2 \le q_3$. We note that $q_1 \ge 0$, since a shortest $(u,v)$-path in $R$ may possibly be the path $uv$ of length~$1$ with no internal vertex. The girth condition implies that $q_1 + q_2 \ge g - 2$ and $q_2 + q_3 \ge g - 2$. We color $u$ and the neighbor of $u$ on each of the paths $Q_2$ and $Q_3$, along with all of $V(G) \setminus V(R)$. Let $Q$ denote the resulting colored set of vertices. Necessarily, $Q$ is a connected forcing set of $G$.

Since every internal vertex of $Q_1$ is uncolored, and only one internal vertex from each of the paths $Q_2$ and $Q_3$ is colored, $q_1 + (q_2 - 1) + (q_3 - 1) + 1 \ge q_1 + q_2 + q_3 - 1$ vertices (including $v$) remain uncolored.
If $q_1 = 0$, then $q_3 \ge q_2 \ge g - 2 - q_1 = g - 2$, implying that the number of uncolored vertices is~$q_1 + q_2 + q_3 - 1 \ge 2(g-2) - 1 = 2g - 5 = (g-2) + (g - 3) \ge g - 2$.
%
If $q_1 \ge 1$, then the number of uncolored vertices is~$q_1 + (q_2 + q_3) - 1 \ge 1 + (g-2) - 1 = g - 2$.
In both cases, $F_{c}(G) \le |Q| \le n - g + 2$.~\qed

\medskip
We establish next the existence of a class of graphs with large connected forcing number.

\begin{figure}[htb]
\begin{center}
\begin{tikzpicture}[scale=.8,style=thick,x=1cm,y=1cm]
\def\vr{2.5pt} 
\path (3,4) coordinate (v);
\path (0,0) coordinate (u0);
\path (0,1) coordinate (u1);
\path (0,2) coordinate (u2);
\path (0,3) coordinate (u3);
\path (-0.5,0) coordinate (u00);
\path (0.5,0) coordinate (u02);

\path (2,0) coordinate (v0);
\path (2,1) coordinate (v1);
\path (2,2) coordinate (v2);
\path (2,3) coordinate (v3);
\path (1.5,0) coordinate (v00);
\path (2.5,0) coordinate (v02);

\path (4,0) coordinate (w0);
\path (4,1) coordinate (w1);
\path (4,2) coordinate (w2);
\path (4,3) coordinate (w3);
\path (3.5,0) coordinate (w00);
\path (4.5,0) coordinate (w02);

\path (6,0) coordinate (x0);
\path (6,1) coordinate (x1);
\path (6,2) coordinate (x2);
\path (6,3) coordinate (x3);
\path (5.5,0) coordinate (x00);
\path (6.5,0) coordinate (x02);

\draw (u00) -- (u1);
\draw (u02) -- (u1);
\draw (u0) -- (u1);
\draw (u1) -- (u2);
\draw (u2) -- (u3);
\draw (u3) -- (v);
\draw (v00) -- (v1);
\draw (v02) -- (v1);
\draw (v0) -- (v1);
\draw (v1) -- (v2);
\draw (v2) -- (v3);
\draw (v3) -- (v);
\draw (w00) -- (w1);
\draw (w02) -- (w1);
\draw (w0) -- (w1);
\draw (w1) -- (w2);
\draw (w2) -- (w3);
\draw (w3) -- (v);
\draw (x00) -- (x1);
\draw (x02) -- (x1);
\draw (x0) -- (x1);
\draw (x1) -- (x2);
\draw (x2) -- (x3);
\draw (x3) -- (v);

\draw (u00) [fill=black] circle (\vr);
\draw (u02) [fill=black] circle (\vr);
\draw (u0) [fill=white] circle (\vr);
\draw (u1) [fill=black] circle (\vr);
\draw (u2) [fill=black] circle (\vr);
\draw (u3) [fill=black] circle (\vr);
\draw (v00) [fill=black] circle (\vr);
\draw (v02) [fill=black] circle (\vr);
\draw (v0) [fill=white] circle (\vr);
\draw (v1) [fill=black] circle (\vr);
\draw (v2) [fill=black] circle (\vr);
\draw (v3) [fill=black] circle (\vr);
\draw (v) [fill=black] circle (\vr);
\draw (w00) [fill=black] circle (\vr);
\draw (w02) [fill=black] circle (\vr);
\draw (w0) [fill=white] circle (\vr);
\draw (w1) [fill=black] circle (\vr);
\draw (w2) [fill=black] circle (\vr);
\draw (w3) [fill=black] circle (\vr);
\draw (x00) [fill=black] circle (\vr);
\draw (x02) [fill=black] circle (\vr);
\draw (x0) [fill=white] circle (\vr);
\draw (x1) [fill=black] circle (\vr);
\draw (x2) [fill=black] circle (\vr);
\draw (x3) [fill=black] circle (\vr);
\end{tikzpicture}
\end{center}

\vskip -0.25 cm
\caption{The graph $G_4$} \label{f:G4}
\end{figure}
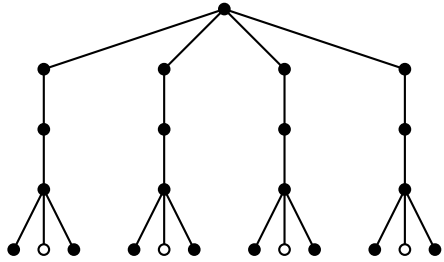

\begin{prop}
For every $\Delta \ge 3$, there exists a connected graph $G_{\Delta}$ of order~$n$ and maximum degree~$\Delta$ such that $F_c(G_{\Delta}) = (\frac{\Delta}{\Delta +1})n + 1$.
\label{p:largeFc}
\end{prop}
\proof Given $\Delta \ge 3$, let $G_{\Delta}$ be obtained from $K_{1,\Delta}$ by subdividing every edge exactly twice, and then attaching $\Delta - 1$ pendant edges to each leaf in the resulting subdivided graph. For example, the graph $G_4$ is illustrated in Figure~\ref{f:G4}, where the darkened vertices form a minimum connected forcing set in $G_4$. The graph $G_{\Delta}$ so constructed has order~$n = \Delta^2 + 2\Delta + 1 = (\Delta + 1)^2$ and maximum degree~$\Delta$. Further, $F_c(G) = n - \Delta = \Delta^2 + \Delta + 1 = (\frac{\Delta}{\Delta +1})n + 1$.~\qed

\medskip

Recall that Amos et al.~\cite{k-Forcing} proved that if $G$ is a general isolate-free graph of order~$n$ with maximum degree~$\Delta$, then $F(G) \le (\frac{\Delta}{\Delta +1})n$. By Proposition~\ref{p:largeFc}, this upper bound fails for $F_c(G)$ when $\Delta \ge 3$.

We establish next the existence of a class of graphs with large connected forcing number in terms of their connected domination number.

\begin{prop}
For all integers $k \ge 1$ and $\Delta \ge 3$, there exists a connected graph $G_{k,\Delta}$ satisfying $F_c(G_{k,\Delta}) = \gamma_c(G)(\Delta - 2) + 2$ and $\gamma_c(G_{k,\Delta}) = k$.
\label{p:conndom}
\end{prop}
\proof Let $k \ge 1$ and $\Delta \ge 3$ be arbitrary given integers. Let $T$ be a tree of order~$k$ satisfying $\Delta(T) < \Delta$. Since $\Delta \ge 3$, we note that such a tree $T$ always exists (as may be seen by taking, for example, $T = P_k$). To each vertex $v$ of $T$, we add $\Delta - d_T(v)$ pendant edges. Let $G$ denote the resulting graph. We note that every vertex in $V(T)$ has degree~$\Delta$ in $G$, while every vertex in $V(G) \setminus V(T)$ is a leaf in $G$. Let $G_{k,\Delta}$ have order~$n$, and note that

\[
\begin{array}{lcl}
n & = & \displaystyle{ |V(T)| + \sum_{v \in V(T)} (\Delta - d_T(v)) } \2 \\
& = & \displaystyle{ k + k \Delta - \sum_{v \in V(T)} d_T(v) } \2 \\
& = & \displaystyle{ k(\Delta + 1) - 2|E(T)| } \2 \\
& = & \displaystyle{ k(\Delta + 1) - 2(k-1) } \2 \\
& = & \displaystyle{ k ( \Delta - 1) + 2 } \2 \\
& = & \displaystyle{ \gamma_c(G_{k,\Delta}) ( \Delta - 1) + 2. }
\end{array}
\]

By construction, if $v \in V(T)$, then $v$ has at least one leaf-neighbor in $G$. Thus, the set $V(T)$ is a minimum connected dominating set of $G$, and so $\gamma_c(G) = |V(T)| = k$.

If $\Delta  = 1$, then $G \cong K_2$, and $F_c(G) = \gamma_c(G) = 1$, implying that $F_c(G) = \gamma_c(G)(\Delta - 2) + 2$. If $\Delta  = 2$, then let $G \cong C_n$, where $n \ge 3$. In this case, $F_c(G) = 2$ and $\gamma_c(G) = n-2$, implying that $F_c(G) = \gamma_c(G) (\Delta - 2) + 2$.

Let $D$ be an arbitrary connected forcing set of $G_{k,\Delta}$.  Necessarily, every vertex in $V(T)$ must be $D$-colored. Further, if some vertex in $V(T)$ has two or more $D$-uncolored leaf-neighbors in $G$, then such a vertex is not a $D$-forcing vertex, implying that $D$ is not a connected forcing set of $G$, a contradiction. Hence, every vertex in $V(T)$ has at most one $D$-uncolored neighbor, implying that at most~$k$ vertices are $D$-uncolored. Therefore, $|D| \ge n - k = (k ( \Delta - 1) + 2) = k = k ( \Delta - 2) + 2 = \gamma_c(G_{k,\Delta}) (\Delta - 2) + 2$. Since $D$ is an arbitrary connected forcing set of $G_{k,\Delta}$, this implies that $F_c(G_{k,\Delta}) \ge \gamma_c(G_{k,\Delta}) (\Delta - 2) + 2$. As shown earlier, $F_c(G) \le \gamma_c(G) (\Delta - 2) + 2$ for all connected graphs $G$ with maximum degree~$\Delta \ge 1$. Consequently, $F_c(G) = \gamma_c(G) (\Delta - 2) + 2$.~\qed

\medskip
We establish next an upper bound on the connected forcing number of a graph in terms of its connected domination number and its maximum degree.

\begin{thm}\label{Forcing Dom}
If $G$ is a connected graph of order~$n$ with maximum degree $\Delta$, then
\[
F_c(G) \le
\left\{
\begin{array}{ll}
\gamma_c(G)(\Delta - 2) + 2 & \mbox{if $\Delta \in \{1,2,n-1\}$ } \2 \\
\gamma_c(G)(\Delta - 1) & \mbox{if $3 \le \Delta \le n-2$.}
\end{array}
\right.
\]
\end{thm}
\proof
Let $G$ be a connected graph of order~$n$ with maximum degree $\Delta \ge 1$. If $\Delta  = 1$, then $G \cong K_2$, and $F_c(G) = \gamma_c(G) = 1$, implying that $F_c(G) = \gamma_c(G)(\Delta - 2) + 2$.
If $\Delta  = 2$, then $G \cong P_n$ or $G \cong C_n$, where $n \ge 3$, and $F_c(G) \le 2$ and $\gamma_c(G) = n-2$, implying that $F_c(G) \le \gamma_c(G) (\Delta - 2) + 2$.
If $\Delta = n-1$, where $n \ge 4$, then $\gamma_c(G) = 1$ and, by Observation~\ref{o:PnCnKn}, $F_c(G) = n-1 = \gamma_c(G) (\Delta - 2) + 2$.
This establishes the bound when $\Delta \in \{1,2,n-1\}$. Hence in what follows, we assume that $3 \le \Delta \le n-2$.

Let $D \subset V$ be a minimum connected dominating set in $G$, and so $\gamma_c(G) = |D|$ and let $\barD = V \setminus D$. Let $m= m(G[D])$ denote the number of edges in $G[D]$. Since $D$ is a connected dominating set, the number of edges in $G[D]$ is at least~$|D| - 1$. We state this formally as follows.

\begin{unnumbered}{Claim~A}
$m \ge |D| - 1$.
\end{unnumbered}

We now color all vertices in the set $D$. Further, for each vertex in $D$, we color all but one neighbor in $\barD$. Let $S$ be the resulting set of colored vertices. Each vertex in $D$ has either all its neighbors $S$-colored or has exactly one $S$-uncolored neighbor and is therefore an $S$-forcing vertex. Since $D$ is a dominating set of $G$, every vertex in $\barD$ has a neighbor in $D$, implying that every vertex of $G$ is $S$-colored or becomes $S$-colored in one forcing step. Thus, the set $S$ is an $S$-forcing set. Suppose that exactly $k$ iterative applications of the $S$-forcing set are required to color all vertices, where each forcing vertex is a vertex of $D$.

There is therefore a subset $D_1$ of $D$ consisting of $|D_1| = k$ vertices and each vertex of $D_1$ is applied in the forcing process. Let $D_2 = D \setminus D_1$, and let $\barD_1$ be the set of vertices in $V \setminus D_1$ that have a neighbor in $D_1$. Further, let $\partial D_2$ be the set of vertices in $D_2$ that have a neighbor in $D_1$. Let $m_i = m(G[D_i])$ for $i \in [2]$, and let $m_{12}$ denote the number of edges between $D_1$ and $D_2$, and so $m_{12} = |[D_1,D_2]| = |[D_1,\partial D_2]|$ and $m = m_1 + m_{12} + m_2$.

We note that the number of $S$-uncolored vertices is at least~$|D_1|$, since each vertex of $D_1$ results in one new vertex changing color in the forcing process. Hence, $F_c(G) \le |S| \le n - |D_1| = |\barD| + |D_2|$.
We state this formally as follows.

\begin{unnumbered}{Claim~B}
$F_c(G) \le |\barD| + |D_2|$.
\end{unnumbered}

Every vertex in $\barD$ is necessarily adjacent to at least one vertex in $D_1$, since after iteratively applying the forcing process using the vertices in $D_1$ all vertices are colored. We proceed further with the following series of claims.

\begin{unnumbered}{Claim~C}
$|\barD| \le |D_1| \cdot \Delta \, - \, 2m_1 \, - \, m_{12}$.
\end{unnumbered}
\proof The number of vertices in $\barD$ is at most the number of edges between $D_1$ and $\barD$. We note that $\barD_1 = \barD \cup \partial D_2$. Thus every vertex in $\barD_1$ belongs to the set $\barD$ or to the set $\partial D_2$, implying that

\[
\begin{array}{lcl}
|\barD| & \le & |[D_1,\barD]| \2 \\
& = & |[D_1,\barD_1]| - |[D_1,\partial D_2]| \2 \\
& = & \displaystyle{ \left( \sum_{v \in D_1} d_{\barD_1}(v) \right) - 2 m_1 - m_{12}} \2 \\
& \le & \displaystyle{ \left( \sum_{v \in D_1} \Delta \right) - 2 m_1 - m_{12} } \2 \\
& = & |D_1| \cdot \Delta \, - \, 2m_1 \, - \, m_{12}. \hspace*{0.25cm} \mbox{{\tiny ($\Box$)}}
\end{array}
\]

\begin{unnumbered}{Claim~D}
$m_{12} + 2m_2 \le |D_2| \cdot \Delta$.
\end{unnumbered}
\proof We note that
\[
|D_2| \cdot \Delta \ge \sum_{v \in D_2} d_{G[D]}(v) =
m_{12} + 2m_2. \hspace*{0.25cm} \mbox{{\tiny ($\Box$)}}
\]

\begin{unnumbered}{Claim~E}
$F_c(G) \le |D| (\Delta - 2) + |D_2| + 2$.
\end{unnumbered}
\proof As observed earlier, $m = m_1 + m_{12} + m_2$. The following now holds by Claim~A, Claim~B, Claim~C and Claim~D.

\[
\begin{array}{lcl}
F_c(G) & \stackrel{(Claim~B)}{\le} & |\barD| + |D_2| \2 \\
& \stackrel{(Claim~C)}{\le} & |D_1| \cdot \Delta - 2m_1 - m_{12} + |D_2| \2 \\
& = & (|D| - |D_2|) \cdot \Delta - 2(m - m_{12} - m_2) - m_{12} + |D_2| \2 \\
& = & (|D| \cdot \Delta - 2m) - |D_2| \cdot \Delta + m_{12} + 2m_2 + |D_2| \2 \\
& \stackrel{(Claim~D)}{\le} & (|D| \cdot \Delta - 2m) + |D_2| \2 \\
& \stackrel{(Claim~A)}{\le} & |D| \cdot \Delta - 2(|D| - 1) + |D_2| \2 \\
& = & |D| (\Delta - 2) + |D_2| + 2.
  \hspace*{0.25cm} \mbox{{\tiny ($\Box$)}}
\end{array}
\]

\begin{unnumbered}{Claim~E}
Every vertex in $D_2$ is a cut-vertex of $G[D]$. \end{unnumbered}
\proof Let $v$ be an arbitrary vertex in $D_2$. Since the set $\barD$ is dominated by the set $D_1$, we note that $D \setminus \{v\}$ is a dominating set of $G$. The minimality of the connected dominating set $D$ therefore implies that the vertex $v$ is necessarily a cut-vertex of $G[D]$.~\smallqed

\begin{unnumbered}{Claim~F}
$|D_1| \ge 2$.
\end{unnumbered}
\proof Suppose, to the contrary, that $|D_1| = 1$. Let $u$ denote the vertex in $D_1$, and so $D_2 = D \setminus \{u\}$. We show that in this case $D_2 = \emptyset$. Suppose that $D_2 \ne \emptyset$.  Let $v$ be a vertex at maximum distance from~$u$ in $G[D]$. By our choice of the vertex $v$, the graph obtained from $G[D]$ by deleting the vertex $v$ is connected, contradicting Claim~E. Therefore, $D_2 = \emptyset$, implying that $\gamma_c(G) = |D| = 1$, a contradiction.~\smallqed

\medskip
We now return to the proof of Theorem~\ref{Forcing Dom}. By Claim~E, $F_c(G) \le |D| (\Delta - 2) + |D_2| + 2$. By Claim~F, $|D_1| \ge 2$, implying that $|D_2|  = |D| - |D_1| \le |D| - 2$, and therefore that $F_c(G) \le |D| (\Delta - 2) + |D| = |D|(\Delta - 1)$.~\qed


\section{Lower Bounds}

In this section we investigate sharp lower bounds on the connected forcing number of a graph. For our first result, recall that a \emph{block} of a graph $G$ is a maximal connected subgraph of $G$ containing no cut vertex of its own. A block of $G$ may, however, contain cut vertices of $G$. Any two blocks of a graph have at most one vertex in common, namely a cut vertex. A block of a graph $G$ containing exactly one cut vertex of $G$ is called an \emph{end block} of $G$. If a connected graph contains a single block, we call the graph itself a \emph{block}. The graph $K_1$ is called the \emph{trivial block}. A \emph{nontrivial block} has order at least~$2$. Every nontrivial block is either $2$-connected or isomorphic to $K_2$. We call a block isomorphic to $K_2$ a $K_2$-\emph{block}. Let $b(G)$ denote the number of $2$-connected blocks of $G$; that is, $b(G)$ is the number of nontrivial blocks of $G$ that are not $K_2$-blocks.

Let $X$ denote the set of cut vertices of a connected graph $G$ and let $Y$ denote the set of its blocks. The \emph{block graph} of $G$ is a bipartite graph $B$ with partite sets $X$ and $Y$ in which a vertex $x \in X$ is adjacent to a vertex $y \in Y$ in $B$ if the block in $G$ corresponding to $y$ contains the vertex corresponding to~$x$.

\begin{thm}\label{Thm:Blocks}
If $G$ is a connected graph, then $F_c(G) \ge b(G) + 1$, and this bound is sharp.
\end{thm}
\proof
If $b(G) = 0$, then $G$ is a tree and the bound follows trivially since in this case $F_c(G) \ge 1 = b(G) + 1$. If $b(G) = 1$, then $G$ is not a tree, and by Observation~\ref{o:path}, $F(G) \ge 2 = b(G) + 1$. Hence, we may assume that $b(G) \ge 2$, for otherwise the desired result follows. Thus, at least two (nontrivial) blocks of $G$ are $2$-connected. Let $G'$ be the graph obtained from $G'$ by iteratively removing all end blocks in $G$ that are $K_2$-blocks. Thus, every end block of $G'$ is $2$-connected and contains exactly one cut vertex of $G$. We note that $b(G') = b(G)$ and $F_c(G) \ge F_c(G')$. Hence it suffices for us to show that $F_c(G') \ge b(G') + 1$.

Let $B'$ be the block graph of $G'$, with partite sets $X$ and $Y$, where $X$ is the set of cut vertices of $G'$ and $Y$ is the set of blocks of $B'$. We observe that the block graph, $B'$, is a tree. Let $Y_1$ be the set of end blocks of $B'$ and let $Y_2$ be the remaining blocks of $B'$, if any. Thus, $|Y| = |Y_1| + |Y_2|$. Counting edges in the tree $B'$, we note that
\[
|X| + |Y_1| + |Y_2| - 1 = |V(B')| - 1 = |E(B')| = |Y_1| + \sum_{v \in Y_2} d_{B'}(v) \ge |Y_1| + 2|Y_2|,
\]
and so, $|X| \ge |Y_2| + 1$. Let $S$ be a minimum connected forcing set of $G'$. We note that every end block of $G'$ has minimum degree at least~$2$, and therefore $S$ must contain at least two vertices from every end block of $G'$ in order for $S$ to be an $S$-forcing set. Further, since every end block of $G'$ is $2$-connected and since $G'[S]$ is connected, the connected forcing set $S$ contains all cut vertices of $G'$. Thus, $X \subset S$ and $S$ contains at least one vertex that does not belong to $X$ from every end block of $G'$. Therefore,
\[
|S| \ge |Y_1| + |X| \ge |Y_1| + |Y_2| + 1 = |Y| + 1 = b(G') + 1.
\]
Thus, $F_c(G) \ge F_c(G') = |S| \ge b(G') + 1 = b(G) + 1$. This establishes the desired upper bound. That the bound is tight may be seen, for example, by taking $k \ge 2$ vertex disjoint cycles (of arbitrary lengths) and identifying one vertex from each cycle into a common vertex $v$ (of degree~$2k$ that belongs to all $k$ cycles). Let $G$ denote the resulting graph. The set $D$ that contains the cut vertex $v$ and a neighbor of $v$ from each of the $k$ cycles is a connected forcing set of $G$, implying that $F_c(G) \le k + 1 = b(G) + 1$. As shown earlier, $F_c(G) \ge b(G) + 1$. Consequently, $F_c(G) \ge b(G) + 1$.~\qed

\medskip
Next we recall a result due to Davila and Kenter \cite{Davila Kenter}, which states that $F(G) \ge \delta + 1$, for graphs with girth $g \ge  4$ and minimum degree $\delta \ge 3$. Since $F_c(G)$ is bounded from below by $F(G)$, as an immediate consequence of this Davila-Kenter  result we have the following lower bound on $F_c(G)$.

\begin{ob} \label{Davila Kenter}
If $G$ is a connected graph with girth $g\ge 4$ and minimum degree $\delta \ge 3$, then $F_c(G) \ge \delta + 1$.
\end{ob}

We next prove a result which relates the girth and minimum degree of a graph to its connected forcing number. We remark that this result is similar to the main conjecture presented in \cite{Davila Kenter}.

\begin{thm}\label{girth lower bound}
If $G$ is a connected graph with girth $g \ge 3$ and minimum degree $\delta \ge 3$, then $F_c(G) \ge \delta + g-3$, and this bound is sharp.
\end{thm}
\proof
If $g = 3$, then $F_{c}(G) \ge F(G) \ge \delta = \delta + g - 3$. If $g=4$, then, by Observation~\ref{Davila Kenter}, $F_c(G) \ge \delta + 1 = \delta + g - 3$. Hence we may assume that $g \ge 5$, for otherwise the desired result is immediate. Let $G$ have order~$n$ and let $S$ be a minimum connected forcing set of $G$, and so $F_c(G) = |S|$. Since $g \ge 5$, we know $G \ncong K_n$ and that $G \ncong K_{1,n-1}$. Hence, by Observation~\ref{o:trivial2} and Theorem~\ref{t:upperbd}, $F_c(G) \le n - 2$, implying that there are at least two $S$-uncolored vertices.
Let $v \in S$ be a vertex that forces initially in the first time step, and $w$ be the non-colored neighbor of $v$ that becomes colored by $v$. Since $\delta \ge 3$, the vertex $w$ has at least two neighbors other than $v$.

Suppose that $w$ has a neighbor $z$, different from $v$, that is $S$-colored. Since $G[S]$ is connected, there is a $(v,z)$-path in $G[S]$. Let $P_{vz}$ be a shortest $(v,z)$-path in $G[S]$. We note that $P_{vz}$ together with $w$ form a cycle $C$, and hence $z$ cannot be a neighbor of $v$ since $g \ge 5$. Since $v$ is an $S$-forcing vertex, all neighbors of $v$ different from $w$ belong to $S$. Further, the vertex $v$ has exactly one neighbor other than $w$ in $C$ since we assumed $P_{vz}$ to be a shortest $(v,z)$-path in $G[S]$. Let $x$ be the neighbor of $v$ different from $w$ that belongs to $C$. Therefore, $V(C) \setminus \{w\} \subset S$ and $N(v) \setminus\{w,x\} \subseteq S$, implying that
\begin{eqnarray*}
|S| & \ge & |V(C) \setminus \{w\}| + |N(v) \setminus\{w,x\}| \\
&= & (|V(C)| - 1) + (d_G(v) - 2) \\
& \ge & (g - 1) + (\delta - 2) \\
& = & \delta + g - 3.
\end{eqnarray*}
Thus, $F_c(G) = |S| \ge \delta + g - 3$. Hence, we may assume that $v$ is the only neighbor of $w$ that is $S$-colored, for otherwise the desired result follows. More generally, we may assume that the non-colored neighbor of an arbitrary $S$-forcing vertex has all its other neighbors non-colored. Since $\delta \ge 3$, this implies that the first two time steps of the forcing process, the first two forcing vertices both belong to $S$. Let $v'$ be the forcing vertex in the second time step, and let $w'$ be the non-colored neighbor of $v'$ that becomes colored by $v'$. Since $S$ is a (connected) forcing set, we may assume, renaming vertices if necessary, that $w$ and $w'$ are adjacent.

Since $G[S]$ is connected, there is a $(v,v')$-path in $G[S]$. Let $P'$ be a shortest $(v,v')$-path in $G[S]$. We note that $P'$ together with the path $v'w'w'v$ form a cycle $C'$ say, and hence $v$ and $v'$ cannot be neighbors since $g \ge 5$. Further, the girth condition implies that $v$ and $v'$ has at most one common neighbor, and such a neighbor necessarily belongs to the cycle $C$. Let $N_v$ and $N_v'$ be the neighbors of $v$ and $v'$, respectively, not on $C$. We note that $N_v \subset S'$, $N_v' \subset S'$, and $N_v \cap N_v' = \emptyset$. Therefore,
\begin{eqnarray*}
|S| & \ge & |V(C) \setminus \{w,w'\}| + |N_v| + |N_v'|  \\
& = & (|V(C)| - 2) + (d_G(v) - 2) + (d_G(v') - 2) \\
& \ge & (g - 1) + 2(\delta - 2) \\
& = & (\delta + g - 3) + (\delta + 3) \\
& \ge & \delta + g - 3.
\end{eqnarray*}
Thus, $F_c(G) = |S'| \ge \delta + g - 3$. This completes the proof of the lower bound. This bound is trivially sharp for $K_n$.~\qed


\section{Open Problems and Conjectures}

As shown in Proposition \ref{p:largeFc}, the connected forcing number can be larger than known upper bounds on the forcing number. Upper bounds on the forcing number of a graph in terms of its order and maximum degree are known. We pose the following question.

\begin{quest}
Given a graph $G$ with order $n$ and maximum degree $\Delta$, is there a function $f(n, \Delta)$ which bounds $F_c(G)$ from above?
\end{quest}

The following conjecture was first given as a lower bound conjecture on $F(G)$ in \cite{Davila Kenter}. Since $F(G)$ is at most $F_c(G)$, and since the original conjecture remains open, we present the following weaker conjecture.

\begin{conj}
If $G$ is a connected graph with minimum degree $\delta$ and girth $g$, then
\begin{equation*}
F_c(G) \ge \delta +(\delta - 2)(g - 3).
\end{equation*}
\end{conj}

For paths, cycles, complete bipartite, and complete graphs, we observe that $F(G)=F_c(G)$. This motivates the following question.

\begin{quest}
What are the necessary and sufficient conditions for $F(G)=F_c(G)$?
\end{quest}

\medskip
We believe that the bound in Theorem~\ref{Forcing Dom} in the case when $3 \le \Delta \le n-2$ cannot be achieved and pose the following question. 

\begin{quest}
Let $G$ be a connected graph with maximum degree $\Delta$, where $3 \le \Delta \le n-1$. Is it true that \[
F_c(G) \le \gamma_c(G)\left(\frac{\Delta^2 - 3\Delta + 3}{\Delta - 1} \right) + 2\left(\frac{\Delta - 2}{\Delta - 1}\right)?
\]
\label{quest3}
\end{quest}

We remark that in the proof of Theorem~\ref{Forcing Dom}, if $|D| \ge (\Delta - 1) |D_2| + 2$, then the upper bound in Question~\ref{quest3} follows from Claim~E. Further, if the upper bound in Question~\ref{quest3} is correct, then it can be shown to be tight.

\end{document}